\documentclass{article}
\usepackage{amssymb}
\usepackage{graphicx}
\usepackage{amsmath}
\numberwithin{equation}{section}
\newtheorem{theorem}{Theorem}[section]         % Theorems numbered within sections
\title{Physics-Informed Neural Networks for the Relativistic Burgers Equation in the Exterior of a Schwarzschild Black Hole}
\author{Shuyang Xiang
	}

\date{\today}
\begin{document}

\maketitle

\begin{abstract}
We introduce a Physics-Informed Neural Networks(PINN) to solve a relativistic Burgers equation in the exterior domain of a Schwarzschild black hole. 

Our main contribution is a PINN architecture that is able to simulate  shock wave formations in such curved spacetime, by training a shock-aware network block and introducing a Godunov-inspired residuals in the loss function. We validate our method with numerical experiments with different kinds of initial conditions. We show its ability  to reproduce both smooth and discontinuous solutions in the context of general relativity. 
\end{abstract}

\section{Introduction}
Recent developments in machine learning and data analysis has led to results across multiple scientific domains. Among these, one of the most outstanding methods is physics informed neural networks (PINN) which incorporate physics laws and partial differential equations directly in the training of the neural network. 
We are interested in exploring the application of PINNs to model compressible fluid flows in the exterior of a Schwarzschild black hole. 

In particular, we consider the relativistic Burgers equation on a Schwarzschild background reading as 
\begin{equation}
\label{eq:relativistic_burgers}
\partial_t \left(\frac {v }{(1-2M/r)^2}\right) + \partial_r\left(\frac{v^2-1}{2(1-2M/r)}\right)=0, \quad r>2M,
\end{equation}
where the unknown velocity $v \in [-1, 1]$ with the normalized  light speed and $M$ the Schwarzschild black hole mass. We assume spherical symmetry for the flow. It is straightforward to recover the standard Burgers equation when the mass $M$ vanishes.  

Previous studies on relativistic fluid dynamics in curved spacetimes, particularly by LeFloch et al.~\cite{Amorim2008,BenArtzi2009}, and on numerical methods by Glimm et al.~\cite{Glimm1965, GlimmMarshallPlohr1984} and Russo et al.~\cite{Russo2005, Semplice2016}, have primarily focused on developing discretization-based schemes to approximate shock wave solutions for the relativistic Burgers equation and the compressible Euler system. In contrast, our goal is to develop a mesh-free numerical approach using PINNs. 

A key challenge is to take the curved spacetime geometry and the possible formation of shock waves into account when designing the network and the training process.

Our main contribution is the development of a robust PINN architecture that handles well shock discontinuities within the relativistic Burgers equation framework on a Schwarzschild spacetime. 

We organize the article as follows: 

In Section 2, we give a summary of the main results~\cite{leflochetal2019} as background on the relativistic Burgers equation Schwarzschild spacetime, including classical steady state and shock solutions. We also give an overview of Physics-Informed Neural Networks (PINNs)~\cite{raissi2019physics}. 
We introduce our relativistic PINNs Section 3. The model is designed to well capture shock formation by incorporating a shock-aware jump block and a Godunov-inspired residual construction. We also introduce the associated loss function and training strategy in detail. Section 4 presents numerical experiments and we demonstrate the efficiency s of our method in learning different type of solutions on a curved spacetime background. 
Finally, in Section 5, we discuss the pros and cos of our methods and the potential of PINNs as a mesh-free alternative for simulating relativistic fluid flows on complicated  metrics.

%========================================================================================================================================
\section{Background}
\subsection{Burgers Equation on Schwarzschild spacetime}
A particular solution of the relativistic Burgers equation~\eqref{eq:relativistic_burgers}  is the steady state solution given by 
\begin{equation}\label{eq:steady-state}
 \partial_r\left(\frac{v^2-1}{2(1-2M/r)}\right)=0, \quad r>2M,
\end{equation}
and we see that such steady state solution has the form
\begin{equation} 
\label{form-steady}
v= \pm \sqrt{1-K\left(1-\frac{2M}{r}\right)},  \quad r>2M
\end{equation}
for $K>1$. 
It is clear that the solution is smooth on $r$ and $\lim_{r \to 2M} v(r) = \pm 1$, i.e., the steady-state solution approaches the speed of light as \( r \to 2M \).. Moreover, we have 
\begin{itemize}
    \item We can define a steady state solution on $(2M,+\infty)$ when $K\in(0, 1]$. 
    \item The steady state solution is defined only on a bounded interval $\left(2M, \frac{2MK^2}{1-K^2}\right)$ when $K>1$. 
\end{itemize}
\begin{figure}[htbp]
    \centering
    \includegraphics[width=0.5\textwidth]{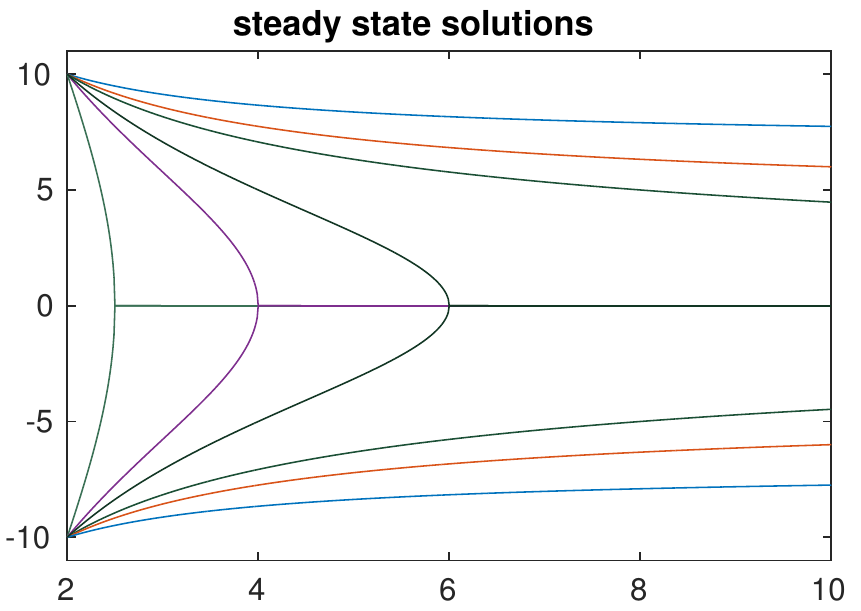}
    \caption{Steady State solutions}
    \label{steadyplot}
\end{figure}
\begin{theorem}[existence theory for the relativistic Burgers model]Consider the
relativistic Burgers equation~\eqref{eq:relativistic_burgers} posed on the outer domain of a Schwarzschild
black hole with mass $M$.  Given a function $v_0: (2M, +\infty) \to (-1, 1)$ with locally finite total variation on $(2M, +\infty)$, that is,
\[
\forall a, b \in (2M, \infty), 
\text{TV}_{[a, b]}(v) < \infty, 
\] then we have: there exists a corresponding
weak solution to $v=v(t, r)$ with locally finite total variation. \end{theorem}

Note that, when talking about the weak solution, we allow the solution to develop discontinuities that satisfy the Rankine–Hugoniot jump condition across shock curves:
\begin{equation}
s[v] =(1-2M/r)\left[\frac{v^2}{2}\right], 
\end{equation}
or, equivalently, $s=(1-2M/r)\frac{v_R+v_L}{2}$ where $s$ is then speed of the shock.  It is obvious that an initial condition with the form of \eqref{form-steady} ensures a unique steady state solution of the same form. In addition to the smooth steady state solutions, we can define a class of
steady shock solutions to the relativistic Burgers equation
\begin{equation}
\label{steady-shock}
v(r)= 
\begin{cases}
\sqrt{1-K\left(1-\frac{2M}{r}\right)}& 2M<r<r_0, \\
-\sqrt{1-K\left(1-\frac{2M}{r}\right)} & r_0>0
\end{cases}
\end{equation} for fixed $r_0>2M$. Since the shock speed $s=0$, the formula \eqref{steady-shock} gives  a weak solution to the relativistic Burgers equation as a steady shock wave. 
\subsection{Physics Informed Neural network}
One of the main challenges in solving partial differential equations (PDEs) is the design of appropriate discretization schemes in me and in space. Poorly designed discretizations can lead to the accumulation of numerical errors at each time step, causing the final numerical solution to diverge. 
Physics-Informed Neural Networks (PINNs) are a class of machine learning models that incorporate physical laws—usually formulated as PDEs—directly into the learning process~\cite{raissi2019physics}. This approach avoids the need for traditional mesh-based discretizations. It gives  promising method to solve complicated PDEs without a complete knowledge of the equation and the associated traditional numerical methods. 
The basic idea is to approximate the solution to a PDE by a neural network, and to define a particular loss function that takes into the equation, the initial condition, boundary condition and any potentially observed data. For a PDE of the form
\begin{equation}
\label{pde}
\mathcal{N}[u](\mathbf{x}, t) = 0,
\end{equation}
we approximate the solution with a neural network $u_\theta(\mathbf{x}, t)$ and the optimal weights  of the network $\theta$ are trained by minimizing the loss function 
\[
\mathcal{L} = \mathcal{L}_{\text{PDE}} + \mathcal{L}_{\text{BC}} + \mathcal{L}_{\text{IC}},
\]
where $\mathcal{L}_{\text{PDE}}$ is the residual of the PDE at given collocation points, and $\mathcal{L}_{\text{BC}}, \mathcal{L}_{\text{IC}}$ are  boundary and initial conditions.
PINNs have been successfully applied to many problems such as fluid mechanics, electromagnetism, etc. 

Extensions and variations have been proposed ever since: Adaptive PINNs (APINNs)  improve convergence by dynamically re-weighting different loss components, handling especially regions with steep gradients~\cite{hu2023apinns}; Extended PINNs(XPINN) applies a domain decomposition and adopt multiple neural networks to solve PDEs in subdomains and helps the case where the geometry is complicated~\cite{jagtap2020xpinns}; Conservative PINNs(cPINNs) introduce the conservation laws into the model and this is of essential importance for hyperbolic systems~\cite{jagtap2021conservative}. 
%=====================================================================================================
\section{Relativistic PINN for shock formation}
We consider the spherical symmetric solution to a relativistic Burgers equation outside of a  Schwarzschild black hole~\eqref{eq:relativistic_burgers} where we recall that $r>2M$ and $M>0$ the black hole mass. 
\subsection{Model Formulation}
We now define a neural network denoted by \( \hat{v}_\theta(t, r) \) that takes as input spacetime coordinates $(t, r) \in [0, \infty)\times (2M, \infty)$ and returns a prediction as the velocity. We want the neural network to minimize the following physics-informed loss:
\begin{equation}
\mathcal{L} = \lambda_{\mathrm{eqn}} \mathcal{L}_{\mathrm{eqn}} + \lambda_{\mathrm{ini}} \mathcal{L}_{\mathrm{ini}} + \lambda_{\mathrm{bnd}} \mathcal{L}_{\mathrm{bnd}},
\end{equation}
where \( \mathcal{L}_{\mathrm{eqn}} \) is the residual of the PDE ~\eqref{eq:relativistic_burgers}, and \( \mathcal{L}_{\mathrm{ini}}, \mathcal{L}_{\mathrm{bnd}} \) are to enforce the initial and boundary conditions.  Note that even the Burgers equation is posed on an unbounded domain, we must truncate the domain to a finite interval $[r_{\text{min}}, r_{\text{max}}]$ which brings the necessity to introduce a boundary condition to anchor the solution in space and ensure numerical stabilities. 

Standard PINN architectures~\cite{raissi2019physics} often fail in the presence of shocks due to the non-smooth nature of the solution. To address this, we propose the following architecture with: 
\begin{itemize}
    \item A standard smooth fully connected network approximating the regular part of the solution (in our case, with 5 hidden layers of the size [16, 32, 32, 32, 16];
    \item A \textbf{shock-aware jump block} that models discontinuities using a dynamically learned shock location.
\end{itemize}

Here, the \texttt{JumpBlock} layer consists of two subcomponents:
\begin{enumerate}
    \item \textbf{Shock location predictor}: a neural network block (with two hidden layers of sizes [32, 32] in our case) that takes time $t$ as input and outputs the position $r_s(t)$ of the shock via a sigmoid activation scaled to the domain:
    \[
    r_s(t) = r_{\min} + (r_{\max} - r_{\min}) \cdot \sigma_{\text{net}}(t), 
    \] where $[r_{\min}, r_{\max}]$ is the truncated domain. 
    \item \textbf{Shock feature generator}: given input $(t, r)$, we consider a smoothed Heaviside function:
    \[
    h(t, r) = \sigma\left(k (r - r_s(t))\right),
    \]
    where $\sigma$ is the sigmoid function and $k$ is a trainable sharpness parameter which controls the steepness of the sigmoid transition.This acts as a differentiable indicator for the location of the shock.
\end{enumerate}

We then define a small neural network contributing to the shock formation:
\begin{equation}
\label{eq:shock-block}
\mathcal{N}_{\text{shock}}(t, r) = \mathcal{F}_{\text{shock}}(t, r, h(t, r)).
\end{equation}

Finally, the complete model for approximation is the sum of the two components:
\[
v_\theta(t, r) = \mathcal{N}_{\text{smooth}}(t, r) + \mathcal{N}_{\text{shock}}(t, r),
\]
where $\mathcal{N}_{\text{smooth}}$ is the fully connected network, and $\mathcal{N}_{\text{shock}}$ is the custom shock block given by \eqref{eq:shock-block}. In the code, such output is written as 
\begin{verbatim}
outputs = out_smooth + jump_out, 
\end{verbatim}
This design gives the network ability to approximate solutions with sharp gradients and discontinuities without diffusions as classical PINNs.

\subsection{Loss Function and Residual Construction}
For a PDE as~\eqref{pde}, the residual is usually given direct by 
\[
||\mathcal{N}[v_\theta](\mathbf{x}, t) ||, 
\]which, unfortunately, has no awareness of the Rankie-Hugoniot condition for our case. 
To address this issue, we adopt a Godunov-inspired initialization.

Indeed, we define the solution residual using a numerical flux as a first-order Godunov scheme. To simply the notation, we now let $(t, r)$ be a space-time input point, and define the metric term:
\[
\text{metric}(r) = 1 - \frac{2M}{r}.
\]
The residual term is computed using:
\begin{itemize}
    \item A data-driven approximation of $\partial_t v$ using automatic differentiation;
    \item A finite-difference approximation of $\partial_r f(v)$, with $f(v) = \frac{v^2 - 1}{2\cdot\text{metric}(r)}$ computed via Godunov fluxes;
    \item The full residual is thus given by :
  \begin{equation}
  \label{residual}
\mathcal{L}_{\text{loss}}(t, r) := \frac{\partial_t v}{\text{metric}(r)^2} + \partial_r f(v).
\end{equation}

\end{itemize}
Moreover, the residual term \eqref{residual} is evaluated via a custom \texttt{GradientLayer} whose inputs are the spacetime coordinates and which returns:
\begin{itemize}
    \item $v(t, r)$: predicted solution;
    \item $\partial_t v(t, r)$: temporal derivative via automatic differentiation;
    \item Left and right interface values $(v_L, v_R)$ via slope limiters ($10^{-5}$ in our application);
    \item $\partial_r f(v)$: computed using the Godunov flux formula:
    \[
    F^G(v_L, v_R) = 
    \begin{cases}
    \max\{f(v_L), f(v_R)\} & \text{if } v_L > v_R, \ s(v_L, v_R) \geq 0, \\
    \min\{f(v_L), f(v_R)\} & \text{if } v_L > v_R, \ s(v_L, v_R) < 0, \\
    f(v_L) & \text{if } v_L \leq v_R, \ v_L \geq 0, \\
    f(v_R) & \text{if } v_L \leq v_R, \ v_R \leq 0, \\
    0 & \text{otherwise},
    \end{cases}
    \]
  where $s = \frac{1}{2}(v_L + v_R) \cdot \text{metric}(r)$ is the shock speed.
\end{itemize}
For the other two terms in the total loss, the initial and boundary conditions, we use the standard $L^2$ loss on predicted $v_\theta(0, r)$, and $v_\theta(t, r_{\min, \max} )$. 

Recall that the equation~\eqref{eq:relativistic_burgers} is defined on \( r \in (2M, +\infty) \), where \( M \) denotes the mass of a black hole and \( r = 2M \) represents the event horizon. We are  particularly interested in the behavior of the solution near the event horizon. However, we avoid sampling training data at $r=2M$ due to the singularity of the Schwarzschild metric at this point. Instead, we consider
rather $r_{\min} =2M + \varepsilon $ for a small $ \varepsilon>0$.  This ensures the boundary data remains finite and the neural network to be trained without numerical instabilities. 
Therefore, boundary conditions are applied at both ends of the radial domain. 

Now let $\mathcal{D}_{\text{eqn}}, \mathcal{D}_{\text{ini}}, \mathcal{D}_{\text{bnd}}$ be collocation points for PDE, initial, and boundary terms respectively. The total loss is written as:
\begin{align}
\label{total-loss}
\mathcal{L} ={}&
\lambda_{\text{eqn}} \cdot \frac{1}{|\mathcal{D}_{\text{eqn}}|} 
\sum_{(t,r)\in\mathcal{D}_{\text{eqn}}} \left| \mathcal{R}(t,r) \right|^2 \notag\\
&+ \lambda_{\text{ini}} \cdot \frac{1}{|\mathcal{D}_{\text{ini}}|} 
\sum_{(0,r)\in\mathcal{D}_{\text{ini}}} \left| v(0,r) - v_0(r) \right|^2 \notag\\
&+ \lambda_{\text{bnd}} \cdot \frac{1}{|\mathcal{D}_{\text{bnd}}|} 
\sum_{(t,r)\in\mathcal{D}_{\text{bnd}}} \left| v(t,r) - v_{\text{bnd}}(t) \right|^2.
\end{align}
where weights $\lambda_{\text{eqn}}, \lambda_{\text{ini}}, \lambda_{\text{bnd}}$ are to be adjusted for training balance.

\subsection{Optimizer and Training Strategy}

To train the PINN model effectively in the presence of sharp features such as shock waves,  we will use a staged training strategy using e L-BFGS-B optimizer from \texttt{scipy.optimize} integrated into our TensorFlow workflow. 

This optimizer is suitable for our medium-sized neural networks (with 4724 parameters according to our design). 
The optimizer parameters were selected as follows: \texttt{factr} = 10.0 (stopping tolerance), memory size \texttt{m} = 50, and maximum line search steps \texttt{maxls} = 100. The optimization proceeds until convergence or a maximum number (500 in our design) of iterations.

We have adopted a three-phase training strategy to improve the convergence: 
\begin{itemize}
    \item \textbf{Warm Start:} Train the model using only initial and boundary data  by setting $\lambda_{\text{eqn}} = 0$ to help the network converge to a physically meaningful solution manifold.
        \item \textbf{Soft Physics Introduction:} Introduce the PDE loss gradually by setting \(\lambda_{\text{eqn}} \ll 1\), allowing the model to adapt to the dynamics without abrupt gradient changes caused by discontinuities and lose the awareness of initial and boundary conditions.
    \item \textbf{Full Physics Training:} Finally, train with all weights set to 1.0, letting the model to approximate a real solution to \eqref{eq:relativistic_burgers}:
    \[
    \lambda_{\text{eqn}} = \lambda_{\text{ini}} = \lambda_{\text{bnd}} = 1.0.
    \]
\end{itemize}

This three-steps training process will significantly improves the convergence of the model that involves shocks. 
\section{Numerical Experiments}

\begin{figure}[htbp]
    \centering
    \includegraphics[width=0.5\textwidth]{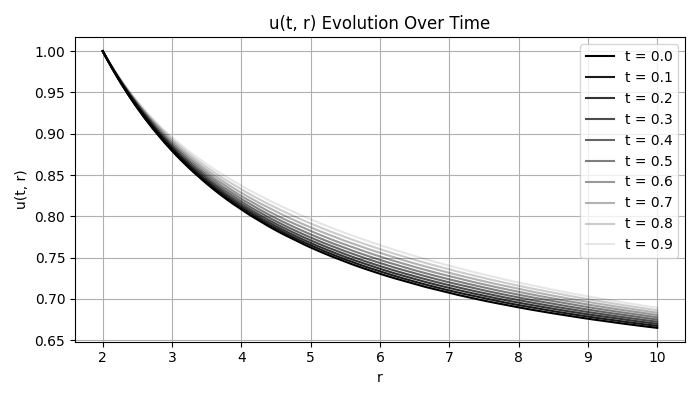}
    \caption{Steady State solution}
    \label{steady-state-figure}
\end{figure}

In this section, we test the ability of our PINN  to model the  dynamics of the relativistic Burgers equation on the Schwarzschild background. We present results in three distinct regimes to highlight the impact of nonlinearity and initial conditions:

We first consider an initially given steady state solution of the form~\eqref{steady-state}. In particular, we solve the equation with \[v_0(r)=\sqrt{1-\frac{3}{4}\left(1-\frac{2M}{r}\right)}\]The model learns the steady state behavior as shown in Figure~\ref{steady-state-figure}. 

\begin{figure}[htbp]
    \centering
    \includegraphics[width=0.5\textwidth]{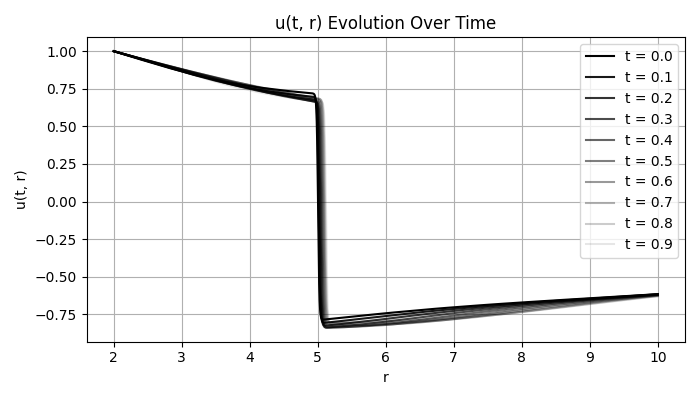}
    \caption{Steady Shock solution}
    \label{steady-shock-figure}
\end{figure}

We next validate the model in a setting of steady shock with the initial condition 
\[
v_0(r)= 
\begin{cases}
\sqrt{1-\frac{3}{4}\left(1-\frac{2M}{r}\right)}& 2M<r<5M, \\
-\sqrt{1-\frac{3}{4}\left(1-\frac{2M}{r}\right)} & r>5M
\end{cases} 
\]The model learns the shock whose shock speed vanishes, see Figure~\ref{steady-shock-figure} . 

\begin{figure}[htbp]
    \centering
    \includegraphics[width=0.5\textwidth]{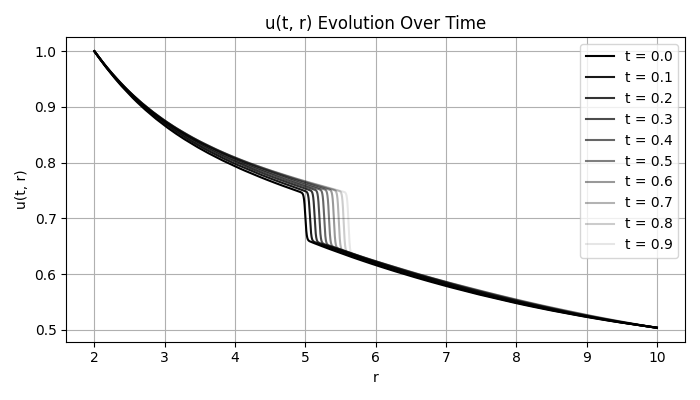}
    \caption{Shock solution}
    \label{shock-figure}
\end{figure}

Finally, we evaluate the model an initial discontinuity which is not steady with \[
v_0(r)= 
\begin{cases}
\sqrt{1-\frac{3}{4}\left(1-\frac{2M}{r}\right)}& 2M<r<5M, \\
\sqrt{1-\frac{15}{16}\left(1-\frac{2M}{r}\right)} & r>5M
\end{cases}
\]. The network successfully captures the shock evolution, as shown in Figure~\ref{shock-figure}. The shock moves to the right hand side gradually. 
These results confirm that the proposed PINN framework, enhanced with an adaptive shock modeling component, is capable of learning both smooth and discontinuous solutions. 

\section {Discussion}
Our relativistic PINN to model compressible fluid flows in the Schwarzschild spacetime shows promising results for different kinds of initial conditions. It successfully integrates the governing Burgers equation directly into the neural network training, enabling the model to capture complex behaviors influenced by curved geometry and shock formation.

While the results are encouraging, we have to admit that the level of accuracy has nothing to do with the  traditional numerical method such as the Glimm method introduced in~\cite{lefloch2018numerical}. However, it still worths mentioning that PINNs offers a possibility of a mesh-free method to the relativistic settings. It is attractive for researchers without  deep background in relativistic fluid dynamics or numerical methods for black hole physics.

The Burgers equation is of course a simple case but the mesh-free nature and the physics-informed learning framework of PINNs gives an alternative to simulate physical systems on a curved geometry background where traditional numerical methods may fail. 

In summary, although PINNs in the context of general relativity currently  stay behind traditional high-accuracy numerical schemes, their flexibility make them an exciting and valuable method and we will continue explore its possibilities. 
\bibliographystyle{plain}
\bibliography{references}

\end{document}